\documentclass[12pt,a4paper]{article}

\usepackage[utf8]{inputenc}
\usepackage[russian]{babel}
\usepackage[left=2cm,right=2cm,
    top=3cm,bottom=3cm,bindingoffset=0cm]{geometry}
    \usepackage{graphicx}
     \usepackage{amsmath,amssymb,euscript,amsthm,amsfonts,mathrsfs,amscd}
\usepackage{color}
\usepackage{amsthm}
\usepackage[affil-it]{authblk}

 \title{Limits of Kalman Filter application in heavy tailed problems}
\author{Valentin Konakov \thanks{Electronic address: \texttt{VKonakov@hse.ru}}}

\author{Pavel Mozgunov%
  \thanks{Electronic address: \texttt{pmozgunov@gmail.com}; corresponding author}}
\affil{International Laboratory of Stochastic Analysis and Its Applications\\ National Research University Higher School of Economics\\ Moscow, Russia}
\date{}
\begin{document}
\maketitle
\begin{abstract}
\footnote{1The article was prepared within the framework of the Academic Fund Program at the National

Research University Higher School of Economics (HSE) in 2014-2015(No 14-05-0007) and supported

within the framework of a subsidy granted to the HSE by the Government of the Russian Federation

for the implementation of the Global Competitiveness Program.} In this paper we consider the behavior of Kalman Filter state estimates in the
case of distribution with heavy tails .The simulated linear state space models with Gaussian measurement noises were used. Gaussian noises in state equation are replaced by components with alpha-stable distribution with di erent parameters alpha and beta. We consider the case when "all parameters are known" and two methods of parameters estimation are compared: the maximum likelihood estimator (MLE) and the expectation- maximization algorithm (EM). It was shown that in cases of large deviation from Gaussian distribution the total error of states estimation rises dramatically. We conjecture that it can be explained by underestimation of the state equation noises covariance matrix that can be taken into account through the EM parameters estimation and ignored in the case of ML estimation.
\end{abstract}
\begin{abstract}
В данной работе рассматривается поведение классического алгоритма Фильтр Калмана в случае, когда ошибки имеют распределение с тяжелыми хвостами. Для этого была использована симулированная модель, в которой шумы наблюдений распределены нормально, а шумы ненаблюдаемых состояний заменены на ошибки, имеющие $\alpha$-устойчивое распределение c двумя меняющимися параметрами $\alpha$ и $\beta$.  Нами рассмотрено два случая: когда все параметры известны, и когда их необходимо оценить. Для оценивания применялись метод максимального правдоподобия и EM-алгоритм. Эмпирически было показано, что в случае больших отклонений от нормальности, средняя ошибка оценивания значений ненаблюдаемых переменных быстро растёт. Мы предполагаем, что это может быть объяснено недооценкой матрицы ковариаций в случае, когда все параметры известны. Но в случае EM-оценивания эта проблема может быть преодолена. Из результатов моделирования следует, что в интервале  $\alpha \in [1.3,2]$ может быть применена стандартная процедура Фильтра Калмана и EM-оценивание. Насколько нам известно, подробный вывод Сглаживания Калмана и EM-алгоритма в русскоязычной литературе ранее не приводился. 
\end{abstract}
\section{Введение}

Рассмотрим дискретную линейную модель пространства-состояний (linear state-space model):
\begin{equation}
 X_{k+1}=A_{k}X_{k}+V_{k+1} \ \ \ \in R^{m}, \ k=0,1,...
 \end{equation}
 \begin{equation}
 Y_{k}=C_{k}X_{k}+W_{k} \ \ \ \ \ \ \ \ \   \in R^{d}, \ \   k=0,1,...
 \end{equation}
где уравнение (1) - уравнение состояний, уравнение (2) - уравнение наблюдений.

\noindent $X_k$ - вектор ненаблюдаемых переменных, $Y_k$ - вектор наблюдаемых переменных.
$A_{k}$ - матрица $m \times m$;  \ \ $V_{k}$ - модельный шум
с функцией плотности $f$; \ \ $C_{k}$ - матрица $d \times m$; \ \  $W_{k}$ - шум наблюдений c функцией плотности $g$; и задано распределение начального вектора $X_0$.

Большое количество моделей (в том числе и со всеми наблюдаемыми переменными) могут быть представлены в виде \textit{линейной} модели пространства состояний, среди них: регресионная модель с меняющимися коэффициентами (TVP), процесс скользящего среднего - MA(q), авторегрессионный процесс - AR(p), модель авторегрессии — скользящего среднего - ARMA (p,q), сезонная модель с шумом, динамическая факторная модель, динамическая факторная модель с общим стохастическим трендом и многие другие модели. 
 
Основной целью при анализе моделей пространства состояний является оценивание функции распределения ненаблюдаемого вектора  $ X_k$ для каждого $k$, основываясь только на зашумленных наблюдениях:
  \begin{equation}
  p(X_k|\mathcal{Y}_{k}),
  \end{equation}
 где  $\mathcal{Y}_{k}$ = $\sigma \{Y_0,...,Y_k \}$ - $\sigma$-алгебра, порожденная наблюдениями $Y_0,...,Y_k$. 

Эта задача была решена Рудольфом Калманом в 1960 году \cite{kalman:1961}, в предположении, что все ошибки и распределение начального вектора гауссовские. Предположение о том, что ошибки в уравнениях гауссовские, существенно, и Фильтр Калмана перестает быть оптимальным в случае не нормальных ошибок.

Целью данной работы является эмпирическая проверка возможности применения Фильтра Калмана для моделей с негауссовскими ошибками. Более точно, было выбрано параметрическое семейство распределений с тяжелыми хвостами, включающее гауссовское распределение, а именно: семейство $\alpha$-устойчивых распределений. В данной работе использованы предварительные результаты, отраженные в \cite{Mozg2014}.

Работа организована следующим образом. Во второй главе будет дан план вывода уравнений Фильтра Калмана и Сглаживания Калмана, затем будет определена альфа-устойчивая случайная величина и описан метод её симулирования, будет введена рассматриваемая в работе модель. В третьей главе дан краткий обзор методов оценивания параметров моделей пространства состояний, часть вывода которых приведена в приложениях. В четвертой главе представлены результаты симуляций. Выводы сделаны в пятой главе. 

\section{Модель}

\subsection{Фильтр Калмана}
Классический Фильтра Калмана (ФК) хорошо известен, тем не менее, для удобства читателя, ниже кратко изложена основная идея данной процедуры.

Рассмотрим уравнения пространства состояний (1) и (2), при следующих предположениях:

$W_0,V_1,W_1,V_2,W_2 \ldots $ - независимые центрированные гауссовские векторы,
$E[V_kV_l^*]=\delta_{kl}Q_k$, $E[W_k W_l^*]=\delta_{kl}R_k$; $ X_0 \sim N(\mu,\Sigma)$ и $E[V_kW_l^*]=0$ $\forall$ $k$ и $l$.

 В этом случае вектора $X_k$,$Y_k$ и $X_k|\mathcal{Y}_k$ для всех моментов времени $k$ будут иметь гауссовское распределение. Так как нормальное распределение полностью определяется двумя первыми моментами, то для оценивания (3) для всех $k$ достаточно найти:
 $$\hat{X}_{k|k}=E[X_{k}|{\mathcal{Y}}_{k}]$$
$${\Sigma}_{k|k}= E[(X_k-\hat{X}_{k|k})(X_k-\hat{X}_{k|k})^*].$$ 
Основная идея Фильтра Калмана (ФК) состоит в разделении процедуры оценивания неизвестного вектора состояний на два этапа: предсказание и корректировка. 
Процедура начинается  с оценивания $\hat{X}_{0|0}$ и ${\Sigma}_{0|0}$, $X_0 \sim N(\mu,\Sigma)$. Используя свойство многомерного нормального распределения, получим:
$$\hat{X}_{0|0}=\mu+\Sigma {C_0}^*(C_0\Sigma {C_0}^* + R_0)^{-1}(Y_0-C_0\mu) $$ 
$${\Sigma}_{0|0}=\Sigma - \Sigma {C_0}^*(C_0\Sigma {C_0}^* + R_0)^{-1}C_0\Sigma,$$
где $(C_0\Sigma {C_0}^* + R_0)^{-1}$ здесь и далее обозначает обычную обратную матрицу, при условии её существования, или, в противном случае, псевдо-обратную.
До поступления наблюдения $Y_{k+1}$, "прогнозируем" значение вектора состояний в момент времени $(k+1)$, используя уравнения $(1)$, то есть:

$$\hat{X}_{k+1|k}=E[X_{k+1}|{\mathcal{Y}}_{k}]=A_k\hat{X}_{k|k}$$
$${\Sigma}_{k+1|k}= E[(A_k(X_k-\hat{X}_{k|k})+V_{k+1})(A_k(X_k-\hat{X}_{k|k})+V_{k+1})^*] = A_k{\Sigma}_{k|k}{A_k}^*+Q_{k+1}$$

Затем поступает наблюдение $Y_{k+1}$.
Обозначим
$$\theta=X_{k+1}-E[X_{k+1}|{\mathcal{Y}}_{k}]= X_{k+1}- \hat{X}_{k+1|k}$$
$$ \xi=Y_{k+1}-E[Y_{k+1}|{\mathcal{Y}}_{k}]= Y_{k+1}- \hat{Y}_{k+1|k}= Y_{k+1}- C_{k+1}\hat{X}_{k+1|k}= \nu_{k+1}$$
Заметим, что $(\theta,\xi)$- гауссовский случайный вектор, независимый от ${\mathcal{Y}}_{k}$.
$$E[\theta|\xi]=E[\theta|\xi,{\mathcal{Y}}_{k}]=E[\theta|\mathcal{Y}_{k+1}]=\hat{X}_{k+1|k+1}- \hat{X}_{k+1|k}$$
Обозначим $\mathbb{V}$ ковариационную матрицу вектора $(\theta,\xi)$:
$$\mathbb{V} =
 \left[ \begin{array}{cc}
 V_{\theta \theta} & V_{\theta\xi}^*\\
 V_{\theta\xi} & V_{\xi\xi}\\
  \end{array} \right]$$
Тогда оценка вектора состояний в момент $(k+1)$ на основе информации $\mathcal{Y}_{k+1}$ запишется следующим образом:
$$ \hat{X}_{k+1|k+1}= \hat{X}_{k+1|k}+E[\theta|\xi]$$
Выражая $E[\theta|\xi]$, как 
$E[\theta|\xi]=E[\theta]+V_{\theta\xi}(V_{\xi\xi})^{-1}(\xi-E[\xi])$, получим:
$$\hat{X}_{k+1|k+1}=\hat{X}_{k+1|k}+G_{k+1}\nu_{k+1},$$ где
$$G_{k+1}=\Sigma_{k+1|k}{C_{k+1}}^*{H_{k+1|k}}^{-1}$$ матрица усиления Калмана (Kalman Gain) и
$$V_{\xi \xi}=E[\nu_{k+1}{\nu_{k+1}}^*]=H_{k+1|k}=C_{k+1}\Sigma_{k+1|k}{C_{k+1}}^*+R_{k+1}$$ 

На заключительном этапе $(k+1)$ шага рассчитывается матрица ошибки оценивания, с учетом полученного наблюдения $Y_{k+1}$:
$$cov(\theta,\theta|\xi)=E[(\theta-E[\theta|\xi]){(\theta-E[\theta|\xi])}^*] = E[(X_{k+1}-\hat{X}_{k+1|k+1}){(X_{k+1}-\hat{X}_{k+1|k+1})}^*] = {\Sigma}_{k+1|k+1}$$
Таким образом, зная оценки на шаге $(k)$, подставляя выражение для $\hat{X}_{k+1|k+1}$ в формулу выше, получаем следующие рекуррентные соотношения для вычисления оценок на $(k+1)$ шаге:

\textbf{Уравнения прогнозирования:}

$\hat{X}_{k+1|k}$=$A_k\hat{X}_{k|k}$

${\Sigma}_{k+1|k}$=$A_k{\Sigma}_{k|k}{A_k}^*+Q_{k+1}$

\textbf{Уравнения корректировки:}

 $\nu_{k+1}$ = $Y_{k+1}- C_{k+1}\hat{X}_{k+1|k}$
 
 $H_{k+1|k}=C_{k+1}\Sigma_{k+1|k}{C_{k+1}}^*+R_{k+1}$

$G_{k+1}=\Sigma_{k+1|k}{C_{k+1}}^*{H_{k+1|k}}^{-1}$

$\hat{X}_{k+1|k+1}=\hat{X}_{k+1|k}+G_{k+1}\nu_{k+1}$

${\Sigma}_{k+1|k+1}$=$(I-G_{k+1}C_{k+1})\Sigma_{k+1|k}$

\subsection{Сглаживание Калмана}
В ходе оценивания ФК наблюдения поступают последовательно, так к шагу $N$ у исследователя имеется наибольший объем информации $ \mathcal{Y}_N$. Чтобы использовать эту информацию введен алгоритм Сглаживания Калмана (СК), который вычисляет: $$\hat{X}_{k|N}=E[X_{k}|{\mathcal{Y}}_{N}]$$
  $${\Sigma}_{k|N}= E[(X_k-\hat{X}_{k|N})(X_k-\hat{X}_{k|N})^*|{\mathcal{Y}}_{N}].$$ 
Понятно, что идея данного алгоритма состоит в получении новых оценок вектора состояний на основе большей информации. Так, если имеется $N$ наблюдений, то оценка ФК для  вектора состояния в момент k+1 вычисляется как условное математическое ожидание относительно $\mathcal{Y}_{k+1}$, а в сглаживании Калмана она вычисляется как условное математическое ожидание относительно $\mathcal{Y}_{N}$, причем $\mathcal{Y}_{N}$ содержит в себе $\mathcal{Y}_{k+1}$.

Вывод формул Сглаживания Калмана приведен в Приложении А.
Сглаженные оценки векторов состояния и сглаженные оценки матрицы ошибок, возникающих при данном оценивании, находятся следующим образом:
$$J_k=\Sigma_{k|k}A_k^*\Sigma_{k+1|k}^{-1}$$
$$\hat{X}_{k|N}=\hat{X}_{k|k}+J_k(\hat{X}_{k+1|N}-\hat{X}_{k+1|k})$$
$$\Sigma_{k|N}=\Sigma_{k|k}+J_k[\Sigma_{k+1|N}-\Sigma_{k+1|k}]J_k^*$$

Как показано в \cite{mozgunov:2014}, доверительный интервал прогноза ненаблюдаемых переменных в ходе оценивания Сглаживанием Калмана \textit{у}же, чем при оценивании Фильтра Калмана, что интуитивно понятно из идеи алгоритма. Более того, ниже будет показано, что в случае гауссовских ошибок СК даёт меньшую среднюю ошибку оценивания, по сравнению с ФК.
\section{Методы оценивания параметров}

Предположим здесь и далее, что параметры модели не зависят от момента времени \textit{k}, тогда вектор параметров, который надо оценить, имеет вид:
$\theta$=
$[\mu,\Sigma,A,C,Q,R]$

\subsection{Метод максимального правдоподобия}

Начнем с некоторой оценки вектора $\theta_0$. Пусть после $j$ шагов имеем оценку $\theta=\theta^{(j)}$ и ищем некоторую лучшую оценку $\hat{\theta}=\theta^{(j+1)}$.
Зная оценку вектора параметров $\theta$ на предыдущем шаге, вычисляем $\nu_k$ и их ковариационные матрицы для $k=1,2, \ldots, N$:
 $$\nu_k(\theta)=y_k-C\hat{X}_{k|k-1},$$ $$H_{k|k-1}(\theta)=C\Sigma_{k|k-1}C^*+R$$
Несложно показать, что $\nu_k$ имеет нормальное распределение со средним $0$ и ковариационной матрицей $H_{k|k-1}$.
Далее строится функция правдоподобия, которая имеет вид: 
\begin{equation}
L_{\nu}(\theta)=\displaystyle\prod_{k=1}^{N}\frac{1}{(2\pi)^{n/2}} |H_{k|k-1}(\theta)|^{-1/2} exp \left(-\frac{1}{2}\nu_k(\theta)^*H_{k|k-1}^{-1}\nu_k(\theta) \right)
\end{equation}

Максимизируя логарифм данного выражения, получаем новую оценку $\hat{\theta}$. Для получения этой оценки часто используют процедуры численной оптимизации, например процедуру Ньютона-Рафсона. Оценивание параметров по методу максимального правдоподобия для моделей пространства состояний производится по следующему алгоритму:
\begin{enumerate}

\item Выбираются некоторые начальные значения для вектора параметров $\theta^{(0)}$.
\item Используя данный вектор параметров, с помощью фильтра Калмана вычисляют $\nu_k(\theta^{(0)})$ и
$H_{k|k-1}(\theta^{(0)})$ для $k=1,2,...,N$. На основе полученных значений вычисляется функция  правдоподобия. 
\item Полученная функция правдоподобия максимизируется (численными методами) по вектору параметров $(\theta)$. В итоге получается некоторый новый вектор параметров $\theta^{(1)}$.
\item Имея новый вектор параметров, повторяются шаги 2-4 до сходимости процедуры.

\end{enumerate}

Данный алгоритм сталкивается с двумя проблемами. Во-первых, из-за численной природы нахождения оценок данный алгоритм весьма времязатратен. Во-вторых, существует риск нахождения локального, а не глобального максимума, поэтому чрезвычайно важен выбор начальных параметров для оценивания. Если же инновации $\nu_k$ имеют не нормальное распределение, то описанный выше алгоритм при некоторых предположениях может давать состоятельные и ассимпотически нормальные оценки, но  уже с другой дисперсией, данный метод называется Quasi MLE.

\subsection{EM-алгоритм}
Впервые данный алгоритм был предложен в \cite{shumway:1982}. Для его вывода введем следующие обозначения.

Пусть на $(\Omega,\mathcal{F},\mathcal{(F}_t),\bar{P})$, где $\bar{P}$ - исходная вероятность(reference probability) случайные  величины $X_k$ и $Y_l$ независимы для любых $k$ и $l$ и подчинены следующим законам распределения.
$$X_k \sim  \mathcal{N}(0,I_n); \ \ Y_l \sim  \mathcal{N}(0,I_m)$$
Рассмотрим новую вероятность $P$ на $(\Omega,\mathcal{F},\mathcal{(F}_t))$, такую что:
$$X_0 \sim N(\mu,\Sigma); \ \ \ X_{k+1}=A_{k}X_{k}+V_{k+1}; \ \ \ V_k \sim N(0,Q_k)$$
$$Y_{k}=C_{k}X_{k}+W_{k}; \ \ \ W_k \sim N(0,R_k)$$
и $\left\lbrace V_k \right\rbrace $ и $\left\lbrace W_l \right\rbrace $ независимы, для всех $k$ и $l$;
Обозначим:
$$ \frac{dP}{d\bar{P}}=\bar{\Lambda}_t,$$где 
$\bar{\Lambda}_t$=$\prod_{k=0}^{t} \bar{\lambda}_k$, (определение $\bar{\lambda}_k$ дано в Приложении Б).
Обозначим:
$$
Q(\theta,\theta')=E_{\theta'}[log\displaystyle\frac{dP_\theta}{dP_{\theta'}}|\mathcal{Y}_N],$$где
вектору  $\theta$=
$[\mu,\Sigma,A,C,Q,R]$ соответствует вероятность  $P_\theta$,
а вектору $\theta'$=
$[\mu',\Sigma',A',C',Q',R']$ соответствует вероятность  $P_\theta'$.

Тогда, обозначая $\bar{E}[\frac{dP_\theta}{d\bar{P}}|\mathcal{Y}_N]=\bar{E}[\bar{\Lambda}_N|\mathcal{Y}_N]$=$L(\theta)$, где $\bar{E}$ - математическое ожидание относительно исходной вероятности $\bar{P}$ и, исходя из неравенства Йенсена, получим следующее неравенство:
$$LogL(\theta)-LogL(\theta')\geq Q(\theta,\theta')$$ 
Основная цель, как и ранее,  состоит в том, чтобы найти такой вектор $\theta$, который максимизирует $L(\theta)$. Начинаем алгоритм с некоторого вектора $\theta^{(0)}$, тогда на шаге $j+1$, $j$ $\geq$0,зная предыдущую оценку $\theta^{(j)}$, ищем вектор $ \theta^{(j+1)}$, так что:
$$\theta^{(j+1)}={\rm arg}\max_{(\theta)} \ Q(\theta,\theta^{(j)})$$
Если удаётся найти такое $\theta^{(j+1)}$, что $Q(\theta^{(j+1)},\theta^{(j)})>0$, то новый вектор увеличивает значение функции правдоподобия по сравнению с предыдущем её значением.
Вывод формул EM-алгоритма представлен в Приложении Б. Формулы для оценивания также могут быть найдены в \cite{shumway:2006}, за исключением формулы для матрицы $C$, так как в \cite{shumway:2006} она предполагается известной, что не всегда так.

Данный алгоритм так же, как и метод максимального правдоподобия, имеет несколько минусов. Во-первых, неизвестны  строгие результаты о сходимости EM-алгоритма. Во-вторых, он очень чувствителен к начальным значениям.
\section{Результаты симуляций}
\subsection{Альфа-устойчивое распределение}
Рассмотрим одномерную $\alpha$-устойчивую величину. Определим её через характеристическую функцию \cite{zolotarev:1986}:
  \begin{equation}
  {\rm log} \phi(t) =  \displaystyle  - \sigma^{\alpha} |t|^{\alpha} \{1-i \beta {\rm sign}(t) {\rm tan} \frac{\pi \alpha}{2} \}+ i \mu t; \alpha \neq 1
  \end{equation}
   \begin{equation}
    {\rm log} \phi(t) =  \displaystyle  - \sigma |t| \{1+i \beta {\rm sign}(t) \frac{2}{\pi} {\rm log}|t| \}+ i \mu t; \alpha = 1
\end{equation}
   где $\alpha \in (0;2], \beta \in [-1,1], \sigma >0, \mu \in R$ и
$\alpha$- характеристическая экспонента $\sigma$- параметр масштаба;
$\beta$- параметр асимметрии. Будем обозначать $\alpha$ - устойчивую случайную величину следующим образом:
$X \sim S_{\alpha}(\sigma,\beta,\mu)$

Класс $\alpha$ - устойчивых распределений был выбран, поскольку случай  $\alpha=2$ и $\beta=0$ ($S_{2}(\sigma,0,\mu))$  соответствует гауссовской случайной величине с параметрами распределения $(\mu,2\sigma^2)$. 

Для симулирования $\alpha$-устойчивой величины мы использовали метод, предложенный в \cite{chambers:1976} и доказанный в \cite{weron:1996}.

Для $\alpha \neq 1$
     \begin{eqnarray}
  X= \displaystyle S_{\alpha , \beta} (\frac{sin \alpha (V + B_{\alpha, \beta})}{(cos V)^{1/\alpha}}) ( \frac{cos(V - \alpha(V + B_{\alpha, \beta}))}{W})^{(1-\alpha)/ \alpha}
\end{eqnarray}
\begin{center}
 $ S_{\alpha , \beta} = \displaystyle [1+\beta^2tan^2 \frac{\pi \alpha}{2}]^{1/(2\alpha)} $\\
  $ B_{\alpha, \beta} = \displaystyle \frac{arctan(\beta tan \frac{\pi \alpha}{2})}{\alpha} $\\
\end{center}
  И для $\alpha = 1$
  \begin{equation}
  X = \displaystyle \frac{2}{\pi} [(\pi /2 + \beta V)tan V - \beta log(\frac{\frac{\pi}{2}W cos V}{(\pi /2) + \beta V})]
\end{equation}
\noindent где V - равномерно распределенная на интервале $[-\pi /2; \pi/2]$ ($V \sim U[-\pi /2; \pi/2])$ случайная величина и W - экспоненциально распределенная случайная величина с параметром 1 ($W \sim Exp (1)$), $V$ и $W$ независимы. 

\subsection{Рассматриваемая модель}
Нами будет симулирована и изучена модель:
 \begin{equation}
 x_{k+1}= Ax_k + \varepsilon_{k+1}
\end{equation} 
 \begin{equation}
 y_{k}= Cx_k + \eta_{k}
\end{equation} 
$\varepsilon_{k}$ $\sim$ $S_{\alpha}(\sigma,\beta,0)$; \ \
$x_0$ $\sim$ $S_{\alpha}(\sigma_2,\beta,\mu)$;\ \
 $\eta_{k}$ $\sim$ $\mathcal{N}(0,R)$, $\varepsilon_{k}$ и $\eta_{l}$ независимы для всех моментов времени $k$ и $l$. 
 
 Мы положим распределение начального вектора состояний $\alpha$-устойчивым с тем же параметром $\alpha$, что и у шума уравнения состояний. Таким образом, благодаря тому, что сумма двух независимых альфа-устойчивых величин с одним и тем же параметром $\alpha$ есть снова альфа-устойчивая величина с тем же параметром $\alpha$, и тому что умножение $\alpha$- устойчивой величины на неслучаную величину изменяет только параметры $\sigma$ и $\mu$, получаем, что вектор состояний во все рассматриваемые моменты времени имеет $\alpha$-устойчивое распределение с одним и тем же параметром $\alpha$. 
 
Ранее мы предполагали, что параметры, фигурирующие в ФК, а именно
$$\theta_k=
[\mu,\Sigma,{A_k},{C_k},{Q_k},{R_k}, S_k]$$
 в каждый момент $k=1,2,\ldots$ известны, но на практике их необходимо оценить. Рассмотрим два метода оценивания параметров рассматриваемого класса моделей.
 
Нами рассмотрено два случая: случай, когда все параметры известны, и случай, когда параметры необходимо оценить.

Для симуляции были выбраны следующие параметры:
 \begin{center}
 $x_{k+1}= x_k + \varepsilon_{k+1}$ ; \ \ $\varepsilon_{k}$ $\sim$ $S_{\alpha}(20,\beta,0)$

$y_{k}= 1.2x_k + \mu_{k}$ ; \ \ $\mu_{k}$ $\sim$ $\mathcal{N}(0,150)$

$x_0$ $\sim$ $S_{\alpha}(50,\beta,100)$
 \end{center}
 Каждая выборка содержит $N=1000$ наблюдений. Для параметров $\alpha$ и $\beta$ было использовано среднее значение $Z$ симуляций ($Z$ будет меняться в зависимости от вычислительных затрат алгоритма).
 \subsection{Все параметры известны}
 Данная глава включает в себя предварительные результаты эмпирического исследования, представленного в \cite{Mozg2014}.
Рассмотрим как изменение $\alpha$ и $\beta$ влияет на среднюю ошибку оценивания ненаблюдаемых состояний. 

 На Рисунке 1 ($Z=10000$) изображена средняя  ошибка оценивания, которая вычисляется следующим образом:
  \begin{equation}
 Error= \displaystyle \frac{1}{Z}  \frac{1}{N} \sum_{m=0}^{Z} \sum_{k=0}^{N}(X^{(m)}_k-\hat{X^{(m)}}_{k|k})^2
   \end{equation}
       \begin{figure}[h!]
 \center\includegraphics[width= 1 \linewidth ]{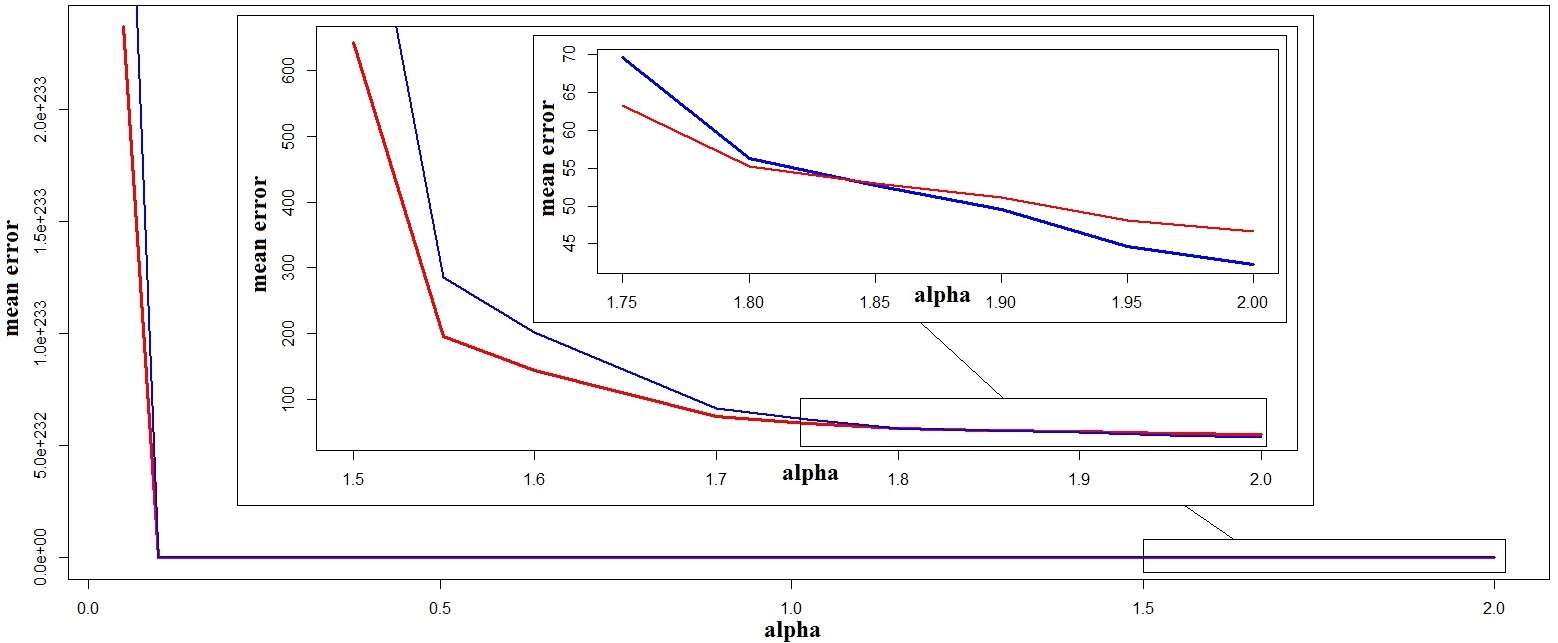}
\caption {Средняя ошибка (Z=10000 симуляций) оценивания ФК (красный) и СК (синий). $\alpha \in [0.05;2]$, шаг=0.05, $\beta = 0$}
\end{figure}

На Рисунке 1 можно видеть, что в случае большого отклонения от гауссовского распределения, средняя ошибка предсказания  растёт, например, для $\alpha = 1.5$ средняя ошибка оценивания в 6 раз больше, чем для гауссовского случая. Из-за огромных значений среднеквадратической ошибки, график средней ошибки в интервале $\alpha \in [0.1:2]$ выглядит как константа, поэтому представлены графики ошибки для меньших интервалов параметра. Заметим, что ФК не так чувствителен к изменению $\alpha$ в интервале $[1.85;2]$, так, например, для $\alpha=1.85$ средняя суммарная ошибка возрастает на 20\%, затем наклон кривой только растёт. Легко видеть, что для этого же диапазона значений $\alpha$ СК даёт меньшую ошибку оценивания по сравнению с ФК. Это объясняется тем, что для нахождения оценок СК, необходимо вновь пользоваться гауссовостью ошибок.

Единственный параметр, который не был установлен истинным - это параметр $\alpha$ распределения ошибок уравнения состояния. Поэтому логично предположить, что большая ошибка связана с несоблюдением предпосылки гауссовости ошибок. Более того, мы предполагаем, что эта ошибка может быть объяснена недооценкой зашумленности ненаблюдаемой переменной. Чтобы это понять, достаточно сравнить $\alpha$ -устойчивое распределение и гауссовское распределение с одним и тем же параметром $\sigma$.\footnote{Гауссовская случайная величина симулировалась методом Бокса-Мюллера, и $\alpha$-устойчивая с.в. симулировалась, как в \cite{weron:1996}}:
  \begin{figure}[h!]
\center\includegraphics[width= 1 \linewidth]{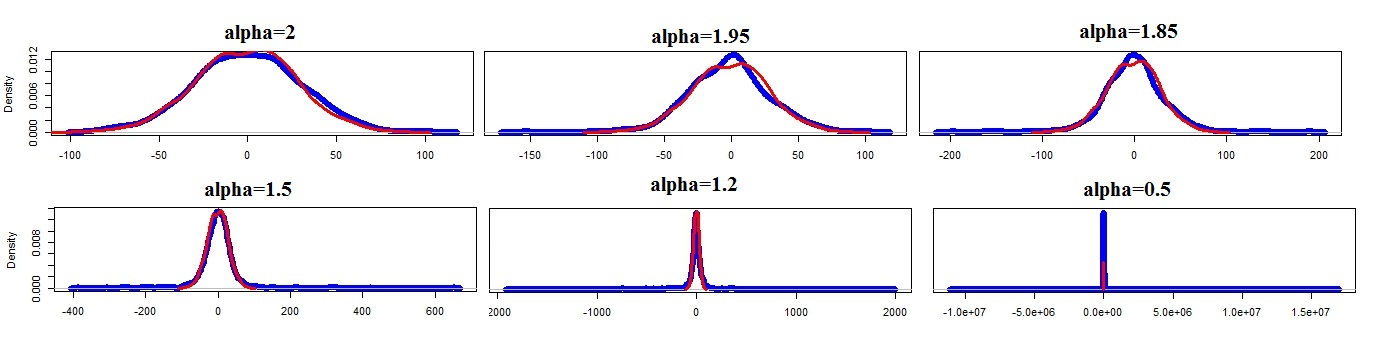}
\caption{Ядерные оценки плотности $\alpha$ - устойчивого распределения (синий)$(\sigma=20;\beta=0)$ и нормального распределения (красный) с тем же параметром $\sigma$.}
\end{figure}
Последствия предположения гауссовости распределения, когда оно на самом деле альфа-устойчивое, представлены на Рисунке 2. Чтобы получить оценки ФК мы заменяем $\alpha$-устойчивое распределение, ядерные оценки плотности которого представлены синим цветом, на гауссовское распределение с таким же параметром $\sigma$, что соответсвует красной линии. Однако, длинные хвосты истинного распределения не покрываются гауссовскими, в результате чего стандартный алгоритм Фильтра Калмана не может распознать прыжки процесса. Легко заметить, что при приближении $\alpha$ к 0, хвосты распределения становятся "тяжелее", и недооценка разброса значений становится всё больше и, как следствие, ошибка возрастает сильнее. В терминах Фильтра Калмана, алгоритм рассматривает уравнение состояний недостаточно зашумленным и приписывает меньший вес (Усиление Калмана) полученному наблюдению.

Средняя общая ошибка для разных параметров $\beta$ и фиксированного $\alpha$ изображена на Рисунке 3.
\begin{figure}[h!]
\center\includegraphics[width=0.8 \linewidth]{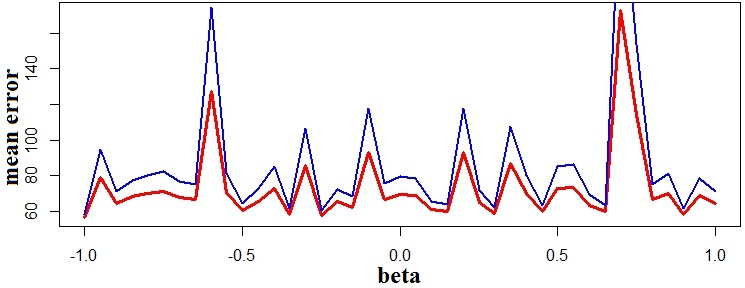}
\caption{Средняя суммарная ошибка (10000 симуляций) оценки  ФК (красный) и СК (синий). $\beta \in [-1;1]$ шаг=$0.05$ $\alpha$ = 1.75}
\end{figure}
 \label{pavel.mozgunov:figure}
  Легко видеть, что случай симметричного распределения ($\beta = 0 $) не соответствует самой маленькой средней ошибке. Сложно сказать, что между ошибкой и параметром $\beta$ есть какая-то зависимость. Более того, рост ошибки не такой большой, как при изменении параметра $\alpha$, поэтому в следующей части мы переключимся только на его изучение.

\subsection{Оценивание параметров}
Рассмотрим поведение двух процедур оценивания, которые были описаны выше: метода максимального правдоподобия (MLE) и EM-алгоритма. Начальные значения взяты на уровне истинных значений, и наша цель изучить, как себя  оценки.\footnote{Из-за вычислительных затрат оптимизации функции правдоподобия, число симуляций было выбрано равным 1000}
 
  \begin{figure}[h!]
\center\includegraphics[width=1 \linewidth]{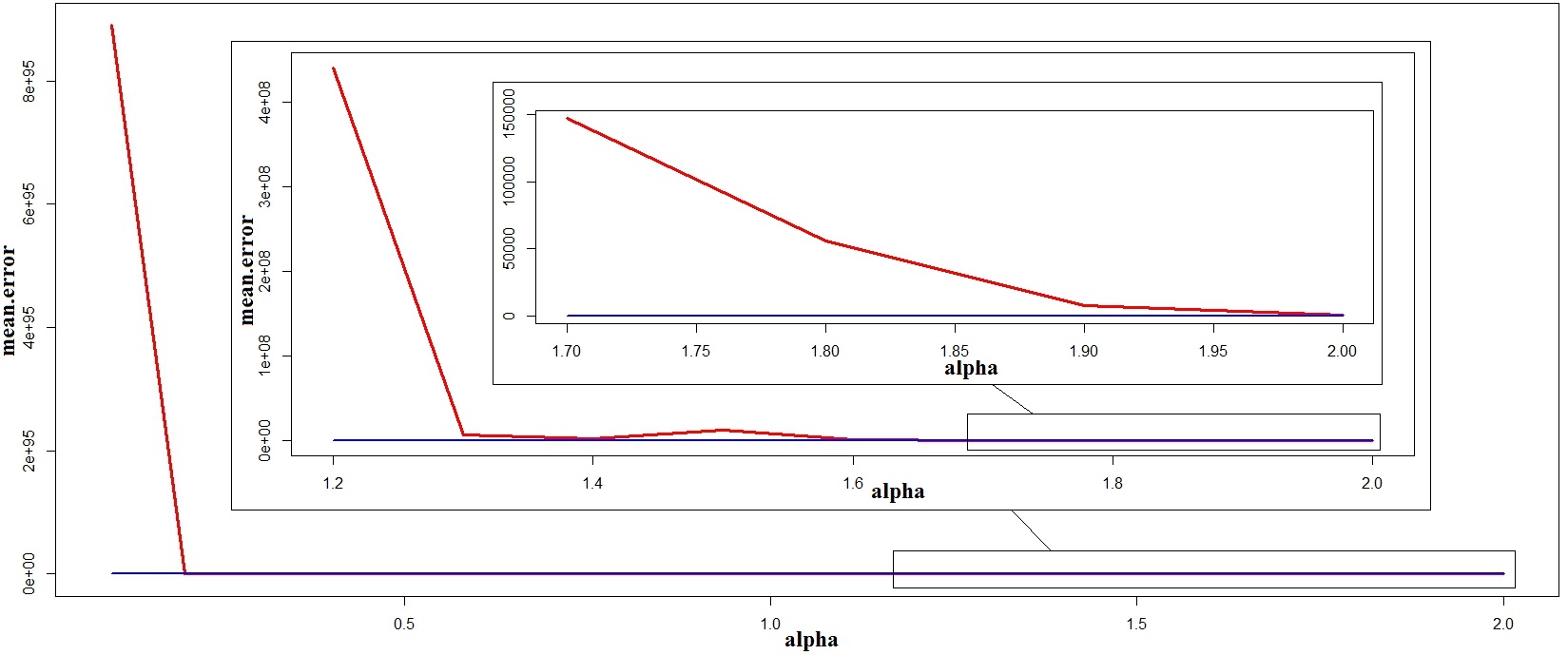}
\caption{Средняя общая ошибка (1000 симуляций).Оценка параметров MLE (красный) и EM (синий).$\alpha \in [0.1;2]$, шаг=0.1 $\beta = 0$}
\end{figure}

На Рисунке 4 можно видеть, что средняя ошибка растёт даже при небольших отклонениях от $\alpha=2$, то есть ML оценивание крайне чувствительно к отклонению от гауссовского случая. Из-за разного порядка ошибки в ходе ML оценивания и EM оценивания, мы изобразили ошибку полученную в ходе оценивания EM-алгоритма на Рисунке 5.

 \begin{figure}[h!]
\center\includegraphics[width=0.77 \linewidth]{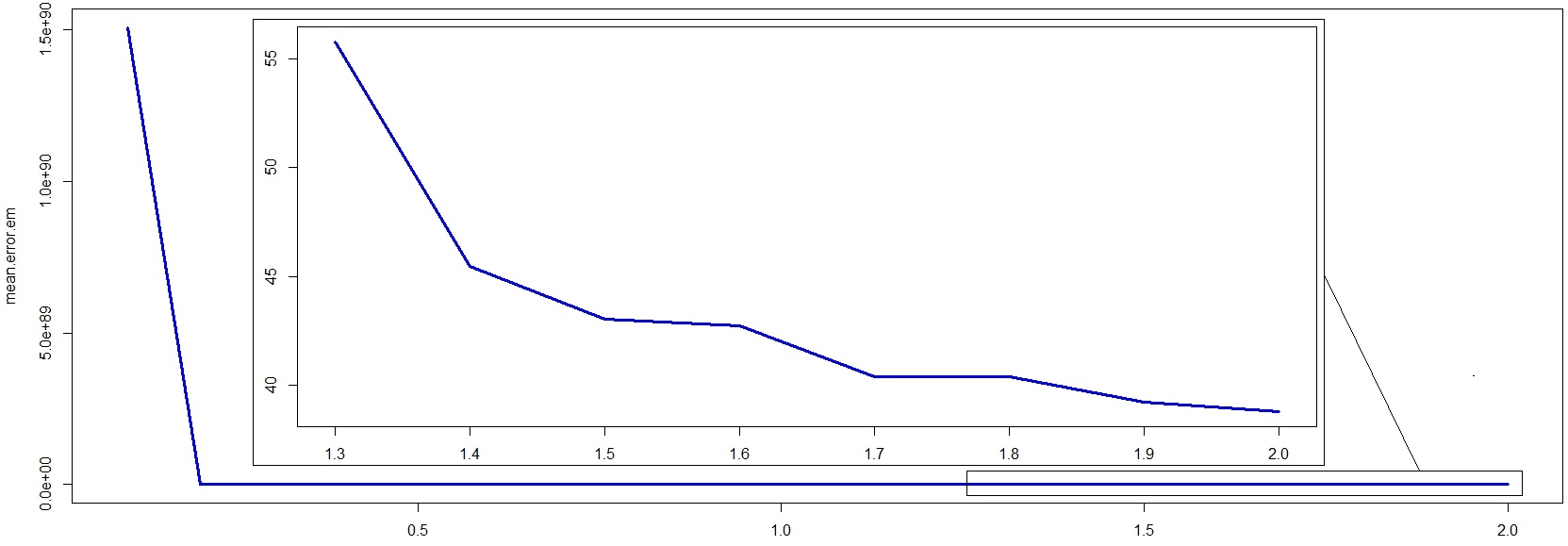}
\caption{Средняя ошибка ФК (1000 симуляции). EM (синий) оценивание параметров. $\alpha \in [0.1;2]$, шаг=0.1 $\beta = 0$}
\end{figure}
 \label{pavel.mozgunov:figure}
 Как видно из Рисунка 5, ошибка оценивания ФК растёт медленно в интервале [1.3;2], например для $\alpha=1.4$ ошибка выросла только  на 12.5\% (в среднем, по сравнению с гауссовским случаем), тогда как без применения алгоритма оценивания она выросла приблизительно в 15 раз, таким образом можно утверждать, что полученый результат уже удовлеторительный. Однако, за пределами данного интервала ошибка уже значительно возрастает, но всё равно меньше, чем ошибка при ML оценивании. Чтобы понять почему EM оценивание даёт более удовлетворительный результат, рассмотрим график оцениваемых параметров (Рисунок 6).

   \begin{figure}[h!]
\center\includegraphics[width=0.5 \linewidth]{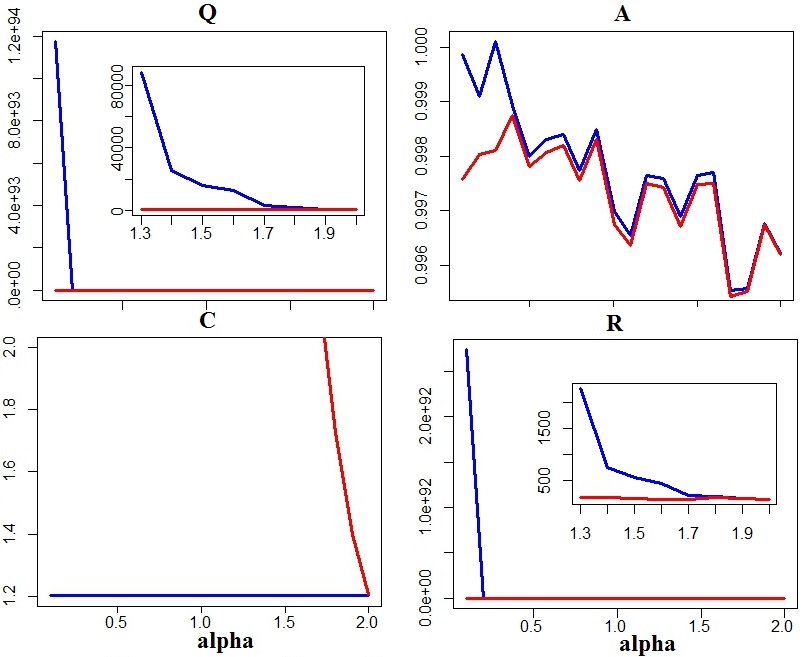}
\caption{ Средние оценки параметров (1000 симуляции) MLE (красный) и EM (синий).$\alpha \in [0.1;2]$, шаг=0.1 $\beta = 0$}
\end{figure}
 \label{pavel.mozgunov:figure}
 
 Как можно видеть из Рисунка 6, EM алгоритм даёт более удовлеторительные результаты оценивания параметров. Алгоритм переоценивает матрицу $Q$ (в сравнении с гауссовским случаем), таким образом, что наблюдения распознаются Фильтром Калмана как более зашумленные, что приводит к большему весу наблюдения в оценивании ненаблюдаемого состояния, что соответсвует увеличению значения параметра $G$ - усиления Калмана.  Аналогичный результат был получен при фиксировании всех остальных параметров и оценивании только матрицы $Q$ с помощью EM алгоритма. Естественно, так как параметр $Q$ теперь больше, это расширяет доверительный интервал нашего прогноза, но в симулированном примере этот рост доверительного интервала выглядит приемлемым, так, например, для  $\alpha=1.4$ мы снизили среднюю ошибку оценивания в более чем в 13 раз, тогда как $\Sigma_{k|k}$ выросла только в $1.15$. К сожалению, данной переоценки недостаточно при больших отклонениях $\alpha$. Однако, мы предъявили достаточно большой интервал возможных значений $\alpha$, для которого можно применить стандартные процедуры оценивания параметров и оценивания  ненаблюдаемых состояний. 

\section{Заключение}
В данной на симулированном примере было показано, что неправильная спецификация шумов уравнения наблюдений значительно увеличивает ошибку оценивания ненаблюдаемых компонент, причем ошибка растёт быстрее чем линейно, при измении параметра $\alpha$. Показано, что EM-алгоритм может оценить дисперсию шума при неверной спецификации ошибок таким образом, что ошибка оценивания ненаблюдаемого вектора возрастает медленее. Так, ошибка оценивания при $\alpha =1.4$ возрастает всего в $1.125$ в сравнении с гауссовским случаем, при приемлемом увеличении доверительного интервала. Данный вывод говорит о том, что для довольно большого интервала значений $\alpha$ может быть использован Фильтр Калмана со стандартной процедурой оценивания параметров EM-алгоритмом.

\section*{Приложение A}
\textbf{Рекурсивная формула для $\hat{X}_{k|N}$.}

Обозначим, $\theta=X_k-\hat{X}_{k|k}$ и $\xi=X_{k+1}-\hat{X}_{k+1|k}$

Заметим, что $\theta$ и $\xi$ образуют гауссовский случайный вектор. Обозначая,
$$\mathbb{V}_1 =
 \left[ \begin{array}{cc}
 V_{\theta \theta} & V_{\theta\xi}^*\\
 V_{\theta\xi} & V_{\xi\xi}\\

  \end{array} \right]$$

$$E[\theta|\xi]=E[\theta]+V_{\theta\xi}{V_{\xi\xi}}^{-1}[\xi-E[\xi]].$$
Здесь
$$V_{\theta\xi} = \Sigma_{k|k} A_k^*$$
$$V_{\xi\xi} = \Sigma_{k+1|k}$$
Тогда
$$X_k-\hat{X}_{k|k}=\Sigma_{k|k}A_k^*\Sigma_{k+1|k}^{-1}(X_{k+1}-\hat{X}_{k+1|k})$$ 

Возьмем условное математическое ожидание $E[...|{\mathcal{Y}}_{N}]$ от левой и правой части и получим:

$$\hat{X}_{k|N}-\hat{X}_{k|k}=J_k(\hat{X}_{k+1|N}-\hat{X}_{k+1|k}),$$ где
$J_k=\Sigma_{k|k}A_k^*\Sigma_{k+1|k}^{-1}$

$$\hat{X}_{k|N}=\hat{X}_{k|k}+J_k(\hat{X}_{k+1|N}-\hat{X}_{k+1|k})=(I-J_kA_k)\hat{X}_{k|k}+J_k\hat{X}_{k+1|N}$$

Так как оценка $\hat{X}_{k|N}$ в выражении выше зависит от $\hat{X}_{k+1|N}$, сглаживание производится начиная с предпоследнего значения. 

\paragraph {Рекурсия для ${\Sigma}_{k|N}$ \\}
Из полученного ранее уравнения:
$$X_k - \hat{X}_{k|N}=X_k-\hat{X}_{k|k}-J_k(X_{k+1|N}-A_k \hat{X}_{k|k} ),$$ имеем:
$$X_k - \hat{X}_{k|N}+J_kX_{k+1|N}=X_k-\hat{X}_{k|k}+J_kA_k \hat{X}_{k|k}.$$ Найдем ковариационные матрицы выражений слева (LHS) и справа(RHS).

$$\displaystyle E((LHS)(LHS)^*) = \Sigma_{k|l} + J_k E(\hat{X}_{k+1|N}\hat{X}_{k+1|N}^*)J_{k}^*.$$
Используя то, что 
$$E(X_{k+1}X_{k+1}^*)= \Sigma_{k+1|l} +  E(\hat{X}_{k+1|N}\hat{X}_{k+1|N}^*)$$
и
$$E(X_{k+1}X_{k+1}^*) = A_kE(X_{k}X_{k}^*)A_k^* + Q_{k+1},$$
получим
$$E((LHS)(LHS)^*) = \Sigma_{k|l} + J_k(A_k E(X_kX_{k}^{*})A_{k}^{*} + Q_{k+1} - \Sigma_{k+1|l})J_{k}^{*}$$
Аналогично,
$$E((RHS)(RHS)^*)= \Sigma_{k|k} + J_kA_k(E(X_kX_{k}^{*}) - \Sigma_{k|k})A_{k}^{*}J_{k}^{*}$$
Приравнивая обе части, получим:
$$\Sigma_{k|l}=\Sigma_{k|k}+J_k[\Sigma_{k+1|l}-\Sigma_{k+1|k}]J_k^*,$$ где $$J_k=\Sigma_{k|k}A_k^*\Sigma_{k+1|k}^{-1}.$$

Заметим, что матрицыы $\Sigma_{k|k}$,$\Sigma_{k+1|k}$ и $J_k$ известны из ФК, и только матрица $\Sigma_{k+1|l}$ неизвестна. То есть в данном случае значение матрицы ошибки в периоде \textit{k}, зависит от её значения в \textit{k+1} периоде, значит нахождение данных сглаженных оценок производится с конца. 

\paragraph {Некоторые полезные результаты из алгоритма сглаживания \\}
Обозначим ковариационную матрицу между парой соседних оценок вектора состояний при сглаживании, следующим образом:

$$\Sigma_{k,l|N}=E[(X_k-\hat{X}_{k|N})(X_l-\hat{X}_{l|N})^*],$$ где $k,l\leq N$

Используя
$$\Sigma_{N,N-1|N}=(I-G_NC_N)A_{N-1}\Sigma_{N-1|N-1}$$ как начальное условие(значение известно из ФК), ковариационную матрицу сглаживания для соседних оценок вектора состояния можно рассчитать как:

$$\Sigma_{k-1,k-2|N}=\Sigma_{k-1|k-1}J_{k-2}^*+J_{k-1}[\Sigma_{k,k-1|N}-A_{k-1}\Sigma_{k-1|k-1}]J_{k-2}^*$$

\section*{Приложение Б}

Обозначим $\mathcal{G}_{k} = \sigma\{X_0,X_1,...,X_k \}$ сигма-алгебру порожденную $X_0,\ldots, X_k$, $\bar{E}$ - математическое ожидание относительно исходной вероятности $\bar{P}$ и $E$ - математическое ожидание относительно вероятности $P$ и
$$\frac{dP}{d\bar{P}}=\bar{\Lambda}_t,$$ где 
$\bar{\Lambda}_t$=$\prod_{k=0}^{t} \bar{\lambda}_k$,  и $\bar{\lambda}_k$, такие что:
\begin{enumerate}
\item $\bar{\lambda}_k$  $\mathcal{G}_{k}$-измерима

\item $\bar{\lambda}_k$>0 c вероятностью 1.

\item $\bar{E}[\bar{\lambda}_{k+1}|\mathcal{G}_{k}] = 1$

\item $\bar{E} [\lambda_0]=1$

\end{enumerate}
Тогда верно следующее соотношение для математического ожидания некоторой случайной $\mathcal{G}_{k}$-измеримой величины A:
$E[A]=\bar{E}[\bar{\Lambda}_tA]$
Пусть $\phi(x)$ гауссовская функция плотности $\mathcal{N}(0,I_n)$, и $\psi(y)$ гауссовская плостность $\mathcal{N}(0,I_m)$. 
Выразим первый член $\bar{\Lambda}_t$, $\bar{\lambda}_0$:
$$\bar{\lambda}_0=\displaystyle\frac{|\Sigma|^{-1/2}\phi(\Sigma^{-1/2}(X_0-\mu))}{\phi(X_0)} \frac{|R_0|^{-1/2}\psi(R_0^{-1/2}(Y_0-C_0X_0))}{\psi(Y_0)} $$
$\bar{\Lambda}_0$=$\bar{\lambda}_0$
Для $k\geq 1$
$$\bar{\lambda}_k=\frac{|Q_k|^{-1/2}\phi(Q_k^{-1/2}(X_k-A_{k-1}X_{k-1})}{\phi(X_k)} \frac{|R_k|^{-1/2}\psi(R_k^{-1/2}(Y_k-C_kX_k))}{\psi(Y_k)} $$
Далее предположим, что параметры модели $\theta_k$ не зависят от момента времени $k$: $\theta_k=\theta$.
\paragraph {E-шаг \\}
$Q(\theta,\theta')=E_{\theta'}[log\displaystyle\frac{dP_\theta}{dP_{\theta'}}|\mathcal{Y}_N]$, где $\theta'$=$\theta^{(j)} $
$$log\displaystyle\frac{dP_\theta}{dP_{\theta'}}=log \bar{\Lambda}_\theta-log\bar{\Lambda}_{\theta'}=\displaystyle\sum_{l=0}^{N}log\bar{\lambda}_l+C(\theta^{\prime}),$$
где $C(\theta^{\prime})$ - константа, не зависящая от вектора параметра $\theta$.

$\displaystyle\sum_{l=0}^{N}log\bar{\lambda}_l + C(\theta^{\prime}) = \displaystyle
-\frac{1}{2}log|\Sigma|-\frac{1}{2}(X_0-\mu)^*\Sigma^{-1}(X_0-\mu)-$

$\displaystyle -\frac{N+1}{2}log|R|$
$\displaystyle -\frac{1}{2}\sum_{k=0}^N(Y_k-CX_k)^*R^{-1}(Y_k-CX_k)-$

$\displaystyle -\frac{N}{2}log|Q|-\frac{1}{2}\sum_{k=1}^N(X_k-AX_{k-1})^*Q^{-1}(X_k-AX_{k-1})+ C(\theta^{\prime})$

Найдем условное математическое ожидание данной величины:

$E_{\theta'}[log\displaystyle\frac{dP_\theta}{dP_{\theta'}}|\mathcal{Y}_N]$=
$\displaystyle -\frac{1}{2}log|\Sigma|-\frac{N+1}{2}log|R|-\frac{N}{2}log|Q| $

$\displaystyle -\frac{1}{2}trace[Q^{-1}E[\sum_{k=1}^N(X_k-AX_{k-1})(X_k-AX_{k-1})^*|\mathcal{Y}_N]] + C(\theta^{\prime})$

$\displaystyle -\frac{1}{2}trace[R^{-1}E[\sum_{k=1}^N(Y_k-CX_{k})(Y_k-CX_{k})^*|\mathcal{Y}_N]]$
$\displaystyle -\frac{1}{2}trace[\Sigma^{-1}E[(X_0-\mu)(X_0-\mu)^*|\mathcal{Y}_N]$

\paragraph {M-шаг \\}

Максимизируем полученную функцию. Перед поиском оптимальных параметров сформулируем две леммы:
\paragraph {\textit{Лемма 1.}}
Если $\bar{Z}=E[Z|\mathcal{Y}]$ и$ V=E[(Z-\bar{Z})(Z-\bar{Z})^*|\mathcal{Y}] $, то справедливо следующее равенство:

$E[ZZ^*|\mathcal{Y}]=V+\bar{Z}\bar{Z}^*$

\paragraph {\textit{Лемма 2.}}
Пусть функция $f(S)=log|S|-trace[SU]$ определена на своей области определения $\mathcal{D}(f)$, состоящей из положительно определенных симметричных матриц $n\times n$, \textit{U} - неотрицательно определенная матрица, тогда $\mathcal{D}(f)$- выпукло, $f$ - вогнута на $\mathcal{D}(f)$и максимум f достигается при $S^{-1}=U$.

Доказательства данных лемм приведены в \cite{shumway:2006}.
\paragraph {Оценка $\mu$ \\}

Согласно Лемме 1:
$E[(X_0-\mu)(X_0-\mu)^*|\mathcal{Y}_N]]$=$\Sigma_{0|N}+(\hat{X}_{0|N}-\mu)(\hat{X}_{0|N}-\mu)^*$,

$trace[\Sigma^{-1}E[(X_0-\mu)(X_0-\mu)^*|\mathcal{Y}_N]=trace[\Sigma^{-1}\Sigma_{0|N}]+
trace[\Sigma^{-1}(\hat{X}_{0|N}-\mu)(\hat{X}_{0|N}-\mu)^*]$.

Так как $\Sigma^{-1}$ и $(\hat{X}_{0|N}-\mu)(\hat{X}_{0|N}-\mu)^*$ две симметричные, неотрицательно определенные матрицы, то справедливо следующее неравенство:
$
trace[\Sigma^{-1}(\hat{X}_{0|N}-\mu)(\hat{X}_{0|N}-\mu)^*]\geq 0.$
Значит, $\hat{\mu}=\hat{X}_{0|N}$.

\paragraph {Оценка матрицы C \\}

$E[(Y_k-CX_{k})(Y_k-CX_{k})^*|\mathcal{Y}_N]$=
$E[Y_kY_k^*|\mathcal{Y}_N]-CE[X_kY_k^*|\mathcal{Y}_N]+CE[X_kX_k^*|\mathcal{Y}_N]C^*-E[Y_kX_k^*|\mathcal{Y}_N]C^* $. 

Заменим:

$\displaystyle E[\sum_{k=0}^NY_kY_k^*|\mathcal{Y}_N] = \textit{K}$;
$\displaystyle E[\sum_{k=0}^NY_kX_k^*|\mathcal{Y}_N] = \textit{-N}$;
$\displaystyle E[\sum_{k=0}^NX_kX_k^*|\mathcal{Y}_N] = \textit{M}$;
$\bar{C}=NM^{-1}$, получим

$E[(Y_k-CX_{k})(Y_k-CX_{k})^*|\mathcal{Y}_N]$=
$(C+\bar{C})M(C+\bar{C})^*+K-\bar{C}M\bar{C}^*$

Далее:

$trace[R^{-1}E[\sum_{k=1}^N(Y_k-CX_{k})(Y_k-CX_{k})^*|\mathcal{Y}_N]]$=
$trace[R^{-1}(C+\bar{C})M(C+\bar{C})^*]+ F$, где $F$ не зависит от $C$.
Cправедливо неравенство:
$trace[R^{-1}(C+\bar{C})M(C+\bar{C})^*] \geq 0 $, то есть $C=-\bar{C}$.
$$\hat{C}=\displaystyle
E[\sum_{k=0}^NY_kX_k^*|\mathcal{Y}_N] \displaystyle(E[\sum_{k=0}^NX_kX_k^*|\mathcal{Y}_N])^{-1}$$

$\displaystyle
E[\sum_{k=0}^NY_kX_k^*|\mathcal{Y}_N]$=$\displaystyle \sum_{k=0}^NE[Y_kX_k^*|\mathcal{Y}_N]$=
$\displaystyle \sum_{k=0}^NY_k\hat{X}_{k|N}^*$

Применяя \textit{Лемму 1}, получим:

$\displaystyle(E[\sum_{k=0}^NX_kX_k^*|\mathcal{Y}_N])$=
$\displaystyle \sum_{k=0}^NE[X_kX_k^*|\mathcal{Y}_N]$=
$\displaystyle \sum_{k=0}^N[\Sigma_{k|N}+\hat{X}_{k|N}{X}_{k|N}^*]$.

Таким образом, в терминах сглаженных оценок Калмана, вычисленных с использованием $\theta$:
$$\hat{C}=\displaystyle \sum_{k=0}^NY_k\hat{X}_{k|N}^*(\sum_{k=0}^N[\Sigma_{k|N}+\hat{X}_{k|N}\hat{X}_{k|N}^*])^{-1}$$

\paragraph {Оценка матрицы A \\}

Для того, чтобы получить оценку теоретического выражения для матрицы $A$, необходимо повторить те же операции, что делались ранее для оценивания матрицы $C$ для следующего выражения:

$E[\sum_{k=1}^N(X_k-AX_{k-1})(X_k-AX_{k-1})^*|\mathcal{Y}_N]$.

В предыдущем алгоритме необходимо заменить $Y_k$ на $X_k$, $C$ на $A$ и $X_k$ на $X_{k-1}$. Получим:
$$\displaystyle \hat{A}=(\sum_{k=1}^N[\Sigma_{k,k-1|N}+\hat{X}_{k|N}\hat{X}_{k-1|N}^*])(\sum_{k=1}^N[\Sigma_{k-1|N}+\hat{X}_{k-1|N}\hat{X}_{k-1|N}^*])^{-1}$$
\paragraph {Оценка матрицы $\Sigma$ \\}

Выпишем часть исходной функции, в которой присутствует искомый параметр $\Sigma$.

$\displaystyle -\frac{1}{2}log|\Sigma|$ $\displaystyle -\frac{1}{2}trace[\Sigma^{-1}E[(X_0-\mu)(X_0-\mu)^*|\mathcal{Y}_N]$.

В обозначениях Леммы 2, $S=\Sigma^{-1}$ и $U=E[(X_0-\mu)(X_0-\mu)^*|\mathcal{Y}_N]$, тогда:

$\displaystyle f(S)=[logS-trace[SU]]$.

По Лемме 2, максимум данной функции достигается при $\Sigma=U$, получим:

$\Sigma=E[(X_0-\mu)(X_0-\mu)^*|\mathcal{Y}_N]$, и по Лемме 1:

$\Sigma=E[(X_0-\mu)(X_0-\mu)^*|\mathcal{Y}_N]=\Sigma_{0|N}+(\hat{X}_{0|N}-\mu)(\hat{X}_{0|N}-\mu)^*$. 
$$\hat{\Sigma}=\Sigma_{0|N}$$
\paragraph {Оценка матрицы R \\}

Выражение части исходной функции, содержащее R, имеет вид:

$\displaystyle -\frac{N+1}{2}log|R|-\frac{1}{2}trace[R^{-1}E[\sum_{k=1}^N(Y_k-CX_{k})(Y_k-CX_{k})^*|\mathcal{Y}_N]]$

В обозначениях Леммы 2, $S=R^{-1}$ и $U=\displaystyle \frac{1}{N+1}E[\sum_{k=1}^N(Y_k-CX_{k})(Y_k-CX_{k})^*|\mathcal{Y}_N]]$

$\displaystyle \frac{N+1}{2}f(S)=\frac{N+1}{2}[log|S|-trace[SU]]$, и максимум функции f(S), достигается при R=U. Подставляя полученное ранее теоретическое выражение для $C$ и беря математическое ожидание, получим:
$$\hat{R}=\displaystyle \frac{1}{N+1}[\sum_{k=0}^NY_kY_k^*-\sum_{k=0}^NY_k\hat{X}_{k|N}^*
(\sum_{k=0}^N[\Sigma_{k|N}+\hat{X}_{k|N}{X}_{k|N}^*])^{-1} \hat{X}_{k|N}Y_k^*]$$

\paragraph {Оценка матрицы Q \\}

Выпишем часть функции, содержащую Q:

$-\frac{N}{2}log|Q|-\frac{1}{2}trace[Q^{-1}E[\sum_{k=1}^N(X_k-AX_{k-1})(X_k-AX_{k-1})^*|\mathcal{Y}_N]]$

Запишем данное выражение в терминах Леммы 2, где $S=Q^{-1}$ и

 $U=\frac{1}{N}E[\sum_{k=1}^N(X_k-AX_{k-1})(X_k-AX_{k-1})^*|\mathcal{Y}_N]]$, тогда оптимальное $Q=U$, то есть:
 
 $$Q=\displaystyle \frac{1}{N}E[\sum_{k=1}^N(X_k-AX_{k-1})(X_k-AX_{k-1})^*|\mathcal{Y}_N]]$$
$Q=\displaystyle
\frac{1}{N}[\sum_{k=0}^N[\Sigma_{k|N}+\hat{X}_{k|N}\hat{X}_{k|N}^*-$

$\displaystyle -
\sum_{k=1}^N
[\Sigma_{k-1,k|N}+\hat{X}_{k-1|N}\hat{X}_{k|N}^*](\sum_{k=0}^N[\Sigma_{k-1|N}+\hat{X}_{k-1|N}\hat{X}_{k-1|N}^*)^{-1}\sum_{k=1}^N
[\Sigma_{k-1,k|N}+\hat{X}_{k|N}\hat{X}_{k-1|N}^*]
]$

\end{document}